\documentclass[11pt]{amsart}
\usepackage{amssymb}
\usepackage{amsmath}

\newcommand{\bbC}{{\mathbb C}}
\newcommand{\bbZ}{{\mathbb Z}}

\newcommand{\bbP}{{\mathbb P}}

\newcommand{\Oh}{{\mathcal O}}

\DeclareMathOperator{\HH}{H}

\DeclareMathOperator{\JJ}{J}

\DeclareMathOperator{\Pic}{Pic}
\DeclareMathOperator{\TT}{T} 
\DeclareMathOperator{\Ext}{Ext}
 
\DeclareMathOperator{\im}{Image}

\DeclareMathOperator{\rank}{rank}

\DeclareMathOperator{\cc}{c}

\DeclareMathOperator{\pf}{pf}

\newcommand{\onto}{\twoheadrightarrow}

\newcommand{\into}{\hookrightarrow}
\newcommand{\by}[1]{\xrightarrow{#1}}

\newcommand{\tensor}{\otimes}
\newcommand{\isom}{\cong}

\newcommand{\sE}{{\mathcal E}}

\newcommand{\sEnd}{\mbox{${\sE}{nd}$}}

\theoremstyle{plain}
\newtheorem{lemma}{Lemma}
\newtheorem{prop}{Proposition}
\newtheorem{thm}{Theorem}

\newtheorem{cor}{Corollary}

\theoremstyle{definition}

\newtheorem{remark}{Remark}

\long\def\comment#1{}
\begin{document}
\title[On some Moduli spaces of stable bundles]{On some Moduli
spaces of stable vector bundles on cubic and quartic threefolds}
\author{Indranil Biswas}
\address{School of Mathematics,
Tata Institute of Fundamental Research,
Homi Bhabha Road, Mumbai 400 005,
INDIA}
\email{indranil@math.tifr.res.in}
\author{Jishnu Biswas}
\address{Theoretical Statistics and Mathematics Unit,
Indian Statistical Institute,
8th Mile Mysore Road,
Bangalore 560059, INDIA}
\email{jishnu@isibang.ac.in}
\author{G. V. Ravindra}
\address{Mathematics Department, Indian Institute of Science,
Bangalore 560 012, INDIA}
\email{ravindra@math.iisc.ernet.in}

\begin{abstract}
We study certain moduli spaces of stable vector bundles of rank two on
cubic and quartic threefolds. In many cases under consideration, it
turns out that the moduli space is complete and irreducible and a
general member has vanishing intermediate cohomology. In one case, all
except one component of the moduli space has such vector bundles. \comment{As a
corollary, : the Fano variety of lines on any smooth cubic or quartic
threefold has the expected dimension.}
\end{abstract}
\maketitle

\section{Introduction}
A vector bundle $E$ on a smooth hypersurface $X\subset\bbP^n$ is
called {\it arithmetically Cohen-Macaulay} (\textit{ACM} for short) if
$$
\HH^i_\ast(X,E):=\bigoplus_{\nu\in\bbZ}\HH^i(X,E(\nu))=0
$$ for all $0<i<n-1$, where $E(\nu):= E\otimes_{{\mathcal O}_X}
   {\mathcal O}_X(\nu)$.  The study of ACM vector bundles of rank two
   on hypersurfaces has received significant attention in the recent
   past (see \cite{Beau, CM1, CM2, CM3, M, MPR1, MPR2, MPR3, R1, BR}.

We note an interesting property of such bundles: Consider the finitely
generated graded $\HH^0_\ast(\bbP^n,\Oh_{\bbP^n})-$module
$\HH^0_\ast(X,E)$.\comment{This is finitely generated.} Any choice of
a minimal set of generators $g_i\in \HH^0(X,E(a_i))$, $1\leq i\leq
  \ell$, gives a minimal resolution
$$
0\to F_1 \by{\Phi} F_0 \to E \to 0\, ,
$$ where $F_0:=\bigoplus_{i=1}^\ell \Oh_{\bbP^n}(-a_i)$ and the map
$F_0\to E$, which yields a surjection $\HH^0_\ast(F_0) \onto
\HH^0_\ast(X,E)$, is defined by the generators $\{g_i\}$. It follows
from the Auslander-Buchsbaum formula (see \cite{E}, Chapter 19) that
$F_1$ is a vector bundle on $\bbP^n$. It is easily verified that $F_1$
is ACM and hence by Horrocks' criterion (see \cite{OSS}, page 39)
splits into a direct sum of line bundles. By changing bases of $F_1$
and $F_0$, the homomorphism $\Phi$ can be chosen to be skew-symmetric,
and it can be checked that the Pfaffian of $\Phi$ is the polynomial
defining $X$ (see \cite{Beau} for details).

We will briefly describe what is known about ACM vector bundles of rank 
two on hypersurfaces. 

$\bullet$ Any ACM vector bundle of rank two on a general hypersurface
$X\subset\bbP^4$ (respectively, $X\subset\bbP^{n\geq 5}$) of degree
$d\geq 6$ (respectively, $d\geq 3$), is a direct sum of line bundles
(see \cite{MPR1, MPR2, R1}). Equivalently, any codimension two {\it
  arithmetically Gorenstein} subscheme (i.e. the zero locus of a
section of an ACM bundle of rank two) in a general smooth hypersurface
$X\subset\bbP^n$ for $n\geq 5$, $d\geq 3$ or $n=4$, $d\geq 6$ is a
complete intersection.

As explained in \cite{R1}, this can be viewed as a ``generalised
Noether-Lefschetz theorem'' for curves in hypersurfaces in
$\bbP^4$. On the other hand, this also provides a verification of the
first non-trivial case of a (strengthening of a) conjecture of
Buchweitz, Greuel and Schreyer (see Conjecture B of \cite{BGS}).

$\bullet$ A general quintic threefold supports only finitely many
indecomposable ACM bundles (see \cite{MPR1}). The zero loci of
non-zero sections of these ACM bundles define non-trivial cycles in
the Griffiths group of the quintic (see \cite{R1}).

$\bullet$ Any ACM bundle on a smooth quadric hypersurface is
isomorphic to a direct sum of line bundles and spinor bundles (see
\cite{knorrer}).

Recall that a vector bundle $E$ over $X$ is called \textit{normalised}
if $h^0(E(-1))=0$ and $h^0(E)\neq 0$. The first Chern class of a
normalised ACM bundle $E$ of rank two satisfies the inequality:
$3-d\leq c_1(E)\leq d-1$, where $d:=\deg{X}$ (see \cite{M,R1}).

Here we take up the study of indecomposable, normalised ACM vector
bundles of rank two on smooth cubic and quartic hypersurfaces in
$\bbP^4$. These hypersurfaces contain a line, which via Serre's
construction (see \cite{OSS}, pages 90--94) yields an indecomposable,
ACM vector bundle of rank two. The condition of being ACM then
implies, by the Grothendieck-Riemann-Roch formula, that the second
Chern class $c_2$ is a function of $c_1$. The main results here are
explicit descriptions of the moduli space of (semi-)stable bundles
with the above Chern classes. This is done by showing that at least in
many cases, any stable bundle of rank two with the above Chern classes
is in fact ACM.

When $c_1(E)=d-1$, the existence of an indecomposable ACM vector
bundle $E$ of rank two on a smooth hypersurface $X$ of degree $d$, is
equivalent to $X$ being the Pfaffian of a skew-symmetric matrix of
size $2d\times 2d$ with linear entries. Starting from this point of
view, the moduli spaces of stable vector bundles with $c_1(E)=d-1$ on
smooth cubic and quartic threefolds have been extensively studied (see
\cite{MTcubic, IMcubic, Druel, IMquartic}). It turns out that in the
cubic case, this moduli space is irreducible and contains only ACM
bundles while in the quartic case there are components of the moduli
spaces which contain only ACM bundles. Using the properties of ACM
bundles, one can compute their dimension and prove smoothness of these
``ACM'' components.

In what follows, we complete the picture by studying ACM bundles for
the remaining (possible) values of $c_1(E)$. Let $X\subset\bbP^4$ be a
smooth hypersurface. The inclusion of $X$ in $\bbP^4$ induces
isomorphisms $\HH^{2i}(X, {\mathbb Z}) = \HH^{2i}(\bbP^4, {\mathbb
  Z})$ for $i=1,2$. Since $\HH^{2i}(\bbP^4, {\mathbb Z})= \mathbb Z$,
where $i=1,2$, both $\HH^{2}(X, {\mathbb Z})$ and $\HH^{4}(X, {\mathbb
  Z})$ will be identified with $\mathbb Z$.  Let
$\mathcal{M}_X(2;c_1,c_2)$ denote the moduli space of normalised,
(semi-)stable vector bundles of rank two with $(c_1(E), c_2(E)) =
(c_1,c_2)$.

The main results are the following:

\begin{thm}\label{cubic}
Let $X\subset\bbP^4$ be a smooth cubic hypersurface.  Let $F(X)$
denote the Fano variety of lines on $X$. There are canonical
isomorphisms
$$\mathcal{M}_X(2;1,2)\isom \mathcal{M}_X(2;0,1)\isom
F(X).$$ Hence these are smooth irreducible surfaces.
\end{thm}

\begin{thm}\label{quartic}
Let $X\subset\bbP^4$ be a smooth quartic hypersurface. Let $F(X)$
denote the Fano variety of lines on $X$ and let $F_2(X)$ denote the
Hilbert scheme of plane conics in $X$.
\begin{enumerate}
\item
Let $X$ be a general quartic threefold.  \comment{Let $\mathcal{Q}$
  denote the Hilbert scheme of all curves in $X$ which are complete
  intersections of three quadrics in $\bbP^4$. The moduli spaces} Let
$\mathcal{M}$ denote any moduli component of $\mathcal{M}_X(2;2,8)$
which contains only ACM bundles. Then $\mathcal{M}$ is smooth of
dimension $5$.\comment{ and the Abel-Jacobi map $c_2:\mathcal{M}\to
  \JJ^2(X)$ is \'etale and quasi-finite.}
\item
There is an isomorphism
$$\mathcal{M}_X(2;1,3)\isom F(X).$$ Hence, for $X$ general, this
moduli space is a smooth irreducible curve.
\item
There is an isomorphism $$\mathcal{M}_X(2;0,2)\isom F_2(X).$$ Hence,
for $X$ general, the moduli space $\mathcal{M}_X(2;0,2)$ is a smooth
irreducible surface.
\end{enumerate}
\end{thm}

\section{The cubic threefold}
\comment{\footnote{We will need to know the following: The Hilbert schemes of
lines, conics, cubics on cubic and quartic hypersurfaces are
irreducible!!. This will imply the irreducibility of the moduli
spaces.}}

Let $X\subset\bbP^4$ be a smooth cubic threefold. 

\begin{lemma}\label{lem1}
Any normalised stable vector bundle $E\to X$ of rank two with
$$
(c_1(E), c_2(E)) = (c_1,c_2)=(1,2)
$$
is indecomposable and ACM.
\end{lemma}

\begin{proof}
Any hypersurface on $X$ is linearly equivalent to $\Oh_X(j)$ for some
$j> 0$. Since $E$ is normalised, we have $\HH^0(X, E(-1))=0$.  Hence
the zero locus of any non-zero section of the rank two vector bundle
$E$ is either empty or pure of codimension $2$. If $E$ has a nowhere
vanishing section, then this implies that $E$ is a direct sum of line
bundles and this contradicts stability.

Let $s\in\HH^0(X,E)$ be a non-zero section whose zero locus $C$ is
pure of codimension $2$ in $X$. So $\deg{C}=c_2(E)=2$. Let $S\subset
X$ be a general hyperplane section so that $\Gamma:=C\cap{S}$ is a
zero-dimensional subscheme of length two. Notice that $\Gamma$ is the
zero locus of the image of $s$ under the map
$\HH^0(E)\into\HH^0(E_{|S})$.

We will show that $E$ is ACM if the restriction $E_{|S}$ is ACM.
To prove this, we proceed as
follows: Consider the short exact sequence
$$
0\to E(n-1) \to E(n) \to E(n)_{|S} \to 0
$$ on $X$.  The homomorphism $\HH^1(X, E(n-1))\longrightarrow
\HH^1(X,E(n))$ in the corresponding long exact sequence of cohomology
groups is a surjection.  Since $\HH^1(X, E(n))=0$ for $n<<0$, this
implies that $\HH^1(X, E(n))=0$ for all $n\in \bbZ$. By Serre duality
we have $\HH^2(X,E(n))= \HH^1(X,E(-3-n))^* = 0$. Hence $E$ is ACM if
$E_{|S}$ is ACM.

To prove that $E_{|S}$ is ACM, we first claim that
$h^2(E_{|S}(-1))=0$. By Serre duality,
$$
h^2(E_{|S}(-1))\,=\, h^0(E_{|S}(-1))\, .
$$
Since $h^0(S, \Oh_S(-1))=0$, in view of the short exact
sequence
\begin{equation}\label{eq1}
0\,\longrightarrow\,\Oh_S\,\stackrel{\cdot s\vert_S}{\longrightarrow}
\, E_{|S} \,\longrightarrow\, I_{\Gamma/S}(1)\,\longrightarrow\, 0
\end{equation}
given by $s\vert_S$, it is enough to show that $h^0(I_{\Gamma/S}) =0$.
\comment{in order to be able to conclude that $h^0(E_{|S}(-1))=0$} Consider
the exact sequence
$$
0\to I_{\Gamma/S} \to \Oh_S \to \Oh_{\Gamma} \to 0\, .
$$ The corresponding induced map $\bbC\isom\HH^0(\Oh_S) \to
\HH^0(\Oh_{\Gamma})\isom\bbC^2$ is injective. Consequently,
$h^2(E_{|S}(-1))=0$ and hence $h^2(E_{|S}(k))=0$ for $k\geq 0$.
\comment{
In this case,
we have $h^0(E_{|S})\geq 3$, so using the short exact sequence
\eqref{eq1}, it follows that $h^0(I_{\Gamma/S}(1))\geq 2$. This in turn
implies that $\Gamma$ is contained in a hyperplane section of $S$ and
so $\Gamma\subset\bbP^2$.}

Now the long exact sequence of
cohomologies associated to the short exact sequence
$$0\to I_{\Gamma/S}(1) \to \Oh_S(1) \to \Oh_{\Gamma}(1) \to 0$$
yields:
$$0\to \HH^0(I_{\Gamma/S}(1)) \to \HH^0(\Oh_S(1)) \to
\HH^0(\Oh_{\Gamma}(1)) \to \cdots .$$ The line bundle $\Oh_S(1)$ being
very ample, separates points and tangents of $S$. Since the length of
the zero-dimensional scheme $\Gamma$ is two, the restriction
homomorphism
$$
\HH^0(S, \Oh_S(1)) \to \HH^0(\Gamma, \Oh_{\Gamma}(1))
$$ is surjective. Thus $h^0(I_{\Gamma/S}(1))= 2$, $h^0(E_{|S})=3$ and
$h^1(I_{\Gamma/S}(1))=0$. Since $h^1(\Oh_S)=0$, it follows from the
short exact sequence \eqref{eq1} that $h^1(E_{|S})=0$. Combining this
with the fact that $h^2(E_{|S}(-1))=0$, we see that $E_{|S}$ is
$1-$regular in the sense of Castelnuovo and Mumford. Hence
$$
h^1(E_{|S}(n))=0
$$
for $n\geq 0$. By Serre duality,
$$h^1(E_{|S}(-n-2))=0  \hspace{3mm} \forall \, n\geq 0.
$$

It remains to be shown that $h^1(E_{|S}(-1))=0$. The Riemann-Roch
Theorem for $E_{|S}$ says
$$\chi(E_{|S})=\frac{3}{2}c_1^2+\frac{3}{2}c_1-c_2+2.$$ Since
$h^0(E_{|S}(-1))=0=h^2(E_{|S}(-1))$ and $c_1(E_{|S}(-1))=-1$,
$c_2(E_{|S}(-1))=2$, by Riemann-Roch we have
$0=\chi(E_{|S}(-1))=h^1(E_{|S}(-1))$. Therefore, $E_{|S}$ is ACM,
completing the proof of the lemma.
\end{proof}

The following lemma gives a minimal resolution of any stable
vector bundle $E$ in the moduli space $\mathcal{M}_X(2;1,2)$.

\begin{lemma}
A vector bundle $E$ in $\mathcal{M}_X(2;1,2)$ has the following
minimal resolution
\begin{equation}\label{minrescubic}
0 \to \Oh_{\bbP^4}(-2)^{\oplus 3}\oplus \Oh_{\bbP^4}(-1) \to
\Oh_{\bbP^4}^{\oplus 3}\oplus \Oh_{\bbP^4}(-1) \to E \to 0.
\end{equation}
\end{lemma}

\begin{proof}
\comment{$E$ is ACM since $I_{\ell/X}$ is ACM. Let $\ell\subset X$ be defined
by $a=b=c=0$ where $a,b,c\in \HH^0(\Oh_X(1))$. One has a diagram}
Take any non-zero section $s$ of $E$.
We have the short exact sequence
\begin{equation}\label{eq2}
0\,\longrightarrow\,\Oh_X\,\stackrel{\cdot s}{\longrightarrow}
\, E \,\longrightarrow\, I_{C/X}(1)\,\longrightarrow\, 0\, ,
\end{equation}
where $C\subset X$ as before is the zero locus of $s$.  Since
$h^0(E)=h^0(E_{|S})=3$, we have $h^0(I_{C/X}(1))=2$ and so
$C\subset\bbP^3$. Now $I_{C/X}$ is ACM $\implies$ $I_{C/\bbP^4}$ is
ACM $\implies$ $I_{C/\bbP^3}$ is ACM. Thus $C\subset\bbP^4$ is a
complete intersection of type $(1,1,2)$ and so there is a surjection
$$
\Oh_X(-1)^{\oplus{2}}\oplus\Oh_X(-2) \onto I_{C/X}\, .
$$

It is straightforward to see this surjection lifts to a map
$$
\Oh_X^{\oplus{2}}\oplus\Oh_X(-1) \rightarrow E\, .
$$
Recall that $E$ has a resolution
$$ 0 \to F_1 \to F_0 \to E \to 0, $$ where $F_0$, $F_1$ are direct
sums of line bundles, and $F_1\isom F_0^\vee(-2)$. Since
$h^0(F_0)=h^0(E)=3$, $h^0(F_1)=0$ and so we have
$F_0=\Oh_{\bbP^4}(-1)^{\oplus{l}}\oplus \Oh_{\bbP^4}^3$ where
$l\in\{1,~3\}$. It is not hard to see that the case $l=3$ is not
possible.
\end{proof}

\begin{lemma}\label{msmooth}
$\HH^2(X,\sEnd{E})=0$.
\end{lemma}

\begin{proof} 
\comment{This follows from results in \cite{MPR1} (see proof of Lemma
\ref{ismooth}). However, we give a direct and elementary proof here.}
We have an exact sequence
$$0\to \Oh_X \to E \to I_{C/X}(1)\to 0.$$ Tensoring this with
$E^\vee$, we get $\HH^2(X,\sEnd{E})\isom \HH^2(E\tensor I_{C/X})$.
{}From the exact sequence
$$0\to E\tensor I_{C/X}\to E \to E_{|C} \to 0$$ we see that
$\HH^2(E\tensor I_{C/X})\isom \HH^1(E_{|C})=\HH^1(N_{C/X})$. By Serre
duality, this is isomorphic to $\HH^0(E(-2)_{|C})$. To compute the
latter, we apply the functor $\mbox{Hom}_{\Oh_{\bbP^4}}(-,\Oh_C(-1))$
to the minimal resolution (\ref{minrescubic}) to get a sequence
$$0 \to \HH^0(E^\vee(-1)_{|C}) \to \HH^0(\Oh_C(-1))^{\oplus{3}} \oplus
\HH^0(\Oh_{C}) \to \HH^0(\Oh_C(1))^{\oplus{3}} \oplus \HH^0(\Oh_{C})\, .
$$ The map $\HH^0(\Oh_C) \to \HH^0(\Oh_C(1))^{\oplus{3}}$ is clearly
injective: this is because the map is induced by linear entries in the
matrix not all of which vanish on $C$. Thus we have
$\HH^0(E^\vee(-1)_{|C})\isom\HH^0(E(-2)_{|C})=0$.
\end{proof}

Let $\mathcal{Q}$ be the Hilbert scheme of all curves in $X$ which are
complete intersections of type $(1,1,2)$ in $\bbP^4$. Consider the
Abel-Jacobi map
$$
\alpha:\mathcal{Q}\to \JJ^2(X)
$$
which associates with the curve $C$ the class of the cycle $3C-2h^2$.
This map $\alpha$ coincides with the composition
$$
\mathcal{Q}\by{p} \mathcal{M}_X(2;1,2) \by{c_2} \JJ^2(X)\, ,
$$
where the map $p$, whose fibre at $[E]\in\mathcal{M}_X(2;1,2)$ is
$\bbP(\HH^0(X,E))\isom\bbP^2$, is given by the Serre construction,
while the map $c_2$ sends any vector bundle $E$ to the image
of the Grothendieck-Chern class $3c_2(E)-2c_1(E)^2$ by the
Abel-Jacobi map.

Since $\HH^1(C,N_{C/X})=0$, the Hilbert scheme
$\mathcal{Q}$ is smooth. Furthermore, by Riemann-Roch, we have
$h^0(N_{C/X})=\deg{N_{C/X}}+2(1-g)$. Here $g=0$ and
$\deg{N_{C/X}}=2$. Thus $h^0(N_{C/X})=4$ and so $\dim(\mathcal{Q})=4$
as a result of which we have $\dim(\mathcal{M}_X(2;1,2))=2$.

\begin{lemma}
The moduli spaces $\mathcal{M}_X(2;1,2)$ and $\mathcal{Q}$ are smooth
and the map
$$
c_2:\mathcal{M}_X(2;1,2)\to \JJ^2(X)
$$
is quasi-finite and \'etale onto its image.
\end{lemma}

\begin{proof}
By Lemma \ref{msmooth},
both $\mathcal{M}_X(2;1,2)$ and $\mathcal{Q}$ are smooth. To
check that the map $c_2$ is smooth, we will use results of Welters
(see \cite{We}): The tangent map
$$
\tau:\TT_{[C]}(\mathcal{Q})\to
\TT_{[3C-2h^2]}\JJ^2(X)
$$
for $\alpha$ is a homomorphism of vector spaces
$\HH^0(C,N_{C/X})\to \HH^1(\Omega^2_X)^\vee$; using Serre duality
and the isomorphism $N_{C/X}^\vee\isom N_{C/X}(-1)$, we get its dual
map
$$
\tau^\ast:\HH^1(\Omega^2_X) \to \HH^1(C,N_{C/X}((-2)).
$$
This homomorphism fits into a commutative diagram:
\[
\begin{array}{ccc}
\HH^0(X,\Oh_X(1)) & \by{R} & \HH^1(\Omega^2_X) \\
{r_C}\downarrow  & & \downarrow{\tau^\ast} \\
\HH^0(C,\Oh_C(1)) & \by{\partial} & \HH^1(C,N_{C/X}(-2)) \\
\end{array}
\]
where $\partial$ is the coboundary map in the long exact
sequence of cohomologies associated to the short exact sequence
$$0\to N_{C/X}\to N_{C/\bbP^4} \to \Oh_C(3) \to 0,$$ and $R$ is the
coboundary map $\HH^0 \to \HH^1$ for the exact sequence
$$0 \to \Omega^2_X \to \Omega^3_{\bbP^4}\tensor N_{X/\bbP^4} \to
\Omega^3_X\tensor N_{X/\bbP^4} \to 0.$$

It is known (see {\it op. ~cit.}) that $R$ is an isomorphism and $r_C$
is a surjection (from a $5$-dimensional vector space onto a
$3$-dimensional vector space). We shall now show that
$\dim{\im(\tau^\ast)}=2$. Since $\im{\tau^\ast}=\im{\partial}$, it is
enough to prove that $\dim{\im(\partial)}=2$. Since \comment{$C$ is
  ACM, this implies that $C$ is a complete intersection in $\bbP^4$
  and hence $N_{C/\bbP^4}$ is a sum of line
  bundles. Since$h^0(I_{C/X}(1))\neq 0$, this implies that
  $C\subset\bbP^3$, hence $N_{C/\bbP^4}$ already has $\Oh_C(1)$ as a
  summand. Now $\deg{C}=c_2(E)=2$, so $C\subset\bbP^3$ has to be a}
$C\subset\bbP^4$ is a complete intersection of type $(1,1,2)$,
$N_{C/\bbP^4}=\Oh_{C}(1)^{\oplus{2}}\oplus\Oh_{C}(2)$. From the
associated long exact sequence of cohomologies one gets a surjection
$\HH^1(C,N_{C/X}(-2)) \onto \HH^1(N_{C/\bbP^4}(-2))\isom
\bbC^2$. Since $\HH^1(N_{C/X}(-2))\isom\HH^0(C,N_{C/X})$ is
$4$-dimensional as seen before, we get $\im(\partial)$ is
$2$-dimensional. Thus the map $\alpha$ has a $2$-dimensional image and
so we are done.
\end{proof}

Let $F(X)$ denote the Fano variety of lines in $X$. It is well known
(see \cite{CG}) that this variety is smooth and irreducible and that
its Albanese is isomorphic to the intermediate Jacobian
$J^2(X)$. Druel \cite{Druel} has proved that the compactified moduli
space $\overline{\mathcal{M}_X(2;2,5)} \isom \widetilde{J^2(X)}$,
where $\widetilde{J^2(X)}$ is the intermediate Jacobian $J^2(X)$ blown
up along a translate of $F(X)$. The following two theorems show that
the remaining two moduli spaces $\mathcal{M}_X(2;1,2)$ and
$\mathcal{M}_X(2;0,1)$ which are of interest to us, also admit an
explicit description thereby proving that they are smooth, and
irreducible.

\begin{thm}\label{firstthm}
The second Chern class map $c_2:\mathcal{M}_X(2;1,2) \to J^2(X)$
induces an isomorphism
$$\mathcal{M}_X(2;1,2)\isom F(X).$$ \comment{Thus the moduli space
$\mathcal{M}_X(2;1,2)$ is irreducible.}
\end{thm}

\begin{proof}
The image of the Abel-Jacobi map $\alpha:\mathcal{Q} \to J^2(X)$ is
known to be isomorphic to $F(X)$. But we have shown above that this
map factors via $\mathcal{M}_X(2;1,2)$. Hence we have the induced map
$c_2:\mathcal{M}_X(2;1,2)\to F(X)$ which is \'etale onto its image. We
give an explicit description of this map and from this deduce that
this is indeed an isomorphism by showing that it has a unique section.

The restriction of the resolution $$0\to F_1\to F_0 \to E \to 0$$
yields a four term exact sequence (see \cite{MPR1}) $$0\to E(-3) \to
\overline{F}_1 \by{\overline{\Phi}} \overline{F}_0 \to E \to 0.$$ The
image $G$ of $\overline\Phi$ is an ACM vector bundle of rank two.  It
can be easily verified that $G(1)$ has a unique section whose zero
locus is a line $\ell$.  Since the minimal resolution of any $E$ as
above is unique up to isomorphism, this map is well defined. To give a
section, we just reverse the process. We first observe that
$\Ext^1_{\Oh_X}(I_{\ell/X}, \Oh_X)\isom \bbC$. This follows by
mimicking the argument given in the proof of Theorem 5.1.1 on pages
94-97 of \cite{OSS}.  As a result, there is a unique vector bundle of
rank two which corresponds to $\ell$ via Serre's construction and
hence the inverse image of $c_2$ has cardinality one. This completes
the proof.
\end{proof}

As a corollary, we obtain
\begin{thm}
$\mathcal{M}_X(2;1,2)\isom \mathcal{M}_X(2;0,1)$.
\end{thm}

\begin{proof}
One checks that $G(1)$ above has Chern classes $(c_1,c_2)=(0,1)$. This
completes the proof.
\end{proof}

\begin{cor}[see also \cite{HRS}, Section $3$]
The Hilbert scheme $\mathcal{Q}$ of plane conics contained in $X$ is
smooth and irreducible of dimension $4$.
\end{cor}

\begin{proof}
The result follows from the fact that $\mathcal{Q}$ is a
$\bbP^2$-bundle over $F(X)$.
\end{proof}

\section{The quartic threefold}
\subsection{Stable Bundles with $(c_1,c_2)=(2,8)$}
The starting point of our investigation is the following result which
can be found in \cite{IMquartic}:

\begin{prop}
Let $X$ be a smooth quartic threefold. Then the moduli space of
normalised, stable vector bundles of rank two with $(c_1,c_2)=(2,8)$
has finitely many components: There is a unique component containing
bundles which via Serre's construction correspond to a union of two
plane sections of $X$, while the remaining are ACM components i.e.,
contain bundles which are ACM (i.e. instantons).
\end{prop}

\begin{lemma}
Let $\mathcal{M}_I$ denote any component of the
moduli space which contains ACM
bundles. Any bundle $E$ which is in $\mathcal{M}_I$ admits a minimal
resolution
\begin{equation}
0 \to \Oh_{\bbP^4}(-2)^{\oplus{4}}\oplus \Oh_{\bbP^4}(-1)^{\oplus{2}}
\to \Oh_{\bbP^4}^{\oplus{4}}\oplus \Oh_{\bbP^4}(-1)^{\oplus{2}} \to E
\to 0.
\end{equation}
\end{lemma}

\begin{proof}
We first show that the zero locus of any non-zero section of $E$ is a
complete intersection of three quadrics in $\bbP^4$. From
Riemann-Roch, it is easy to see that $E$ has $4$ sections, hence $E$
is at least $4$-generated. Let $C\subset X$ be the zero locus of any
non-zero section of $E$. From the short exact sequence
$$ 0\longrightarrow \Oh_X \stackrel{\cdot s}{\longrightarrow} E
\longrightarrow I_{C/X}(2)\longrightarrow 0
$$ we see that $h^0(I_{C/X}(2))=3$, and hence $C\subset X$ is
contained in three mutually independent quadrics $Q_i\subset\bbP^4$
for $1\leq i \leq 3$. Since $\deg{C}=8$ and $C\subset
\cap_{i=1}^3{Q}_i$, it follows that $C=\cap_{i=1}^3{Q}_i$. Next, we
shall show that $E$ is $6$-generated with the other two generators in
degree $1$. To see this notice that since $F_1\isom F_0^\vee(-2)$, any
resolution must be of the form
\begin{equation}
0 \to \Oh_{\bbP^4}(-2)^{\oplus{4}}\oplus\Oh_{\bbP^4}(-1)^{\oplus{k}} \to
\Oh_{\bbP^4}^{\oplus{4}}\oplus \Oh_{\bbP^4}(-1)^{\oplus{k}} \to E \to
0,
\end{equation}
where $k\in\{0,2,4\}$. If $k=0$, then this implies that $E$ is
$0$-regular, and hence $h^0(E)=0$ which is a contradiction.  When
$k=4$, then it can be easily seen that $X$ is not smooth.
\end{proof}

The following result shows that when $X$ is general, the ``ACM''
components of rank two stable bundles that we are interested in, are
all smooth. Note that it is easy to prove this statement for each case
by using the minimal resolutions as in the cubic case (see Lemma
\ref{msmooth}). 

\begin{lemma}\label{ismooth}
Let $X$ be a general quartic threefold. For any ACM rank two bundle
$E$,
$$\HH^2(X,\sEnd{E})=0\, .
$$ \comment{So $\mathcal{M}$ is smooth at the
point $[E]$.}
\end{lemma}

\begin{proof}
This follows from results in \cite{MPR1}. First one shows that
 $\HH^2_\ast(X,\sEnd{E})$ is a finite length module over the graded
 ring $\HH^0_\ast(X,\Oh_X)$ with the generator living in the
 one-dimensional vector space $\HH^2(X,(\sEnd{E})(-d))$ (see
 \cite[Corollary 2.3 and sequence (6)]{MPR1}.  Next it is shown that
 for any $g\in\HH^0(\Oh_{X}(d))$, the multiplication map
 $\HH^2(X,(\sEnd{E})(-d)) \by{g} \HH^2(X,\sEnd{E})$ is zero
 \cite[Corollary 3.8]{MPR1}.
\end{proof}

Let $\mathcal{M}$ denote the union of the ACM components
$\mathcal{M}_I$. By results in section $3$, \cite{MPR1}, there are
only finitely many such components. This proves in particular that
$\mathcal{M}$ is a finite disjoint union of smooth varieties.

\begin{thm}
Let $X$ be a general quartic threefold.  Let $\mathcal{Q}$ denote the
Hilbert scheme of all curves in $X$ which are complete intersections
of three quadrics in $\bbP^4$. The moduli spaces $\mathcal{M}$,
$\mathcal{Q}$ are smooth of dimension $5$ and $8$ respectively and the
map $c_2:\mathcal{M}\to \JJ^2(X)$ is \'etale and quasi-finite.
\end{thm}

\begin{proof}
 When $X$ is general, by Lemma \ref{ismooth}, $\mathcal{M}$ is
  smooth. Since $\mathcal{Q}$ is a projective bundle over
  $\mathcal{M}$, it is also smooth.
  To check that the map $c_2$ is smooth, we will again use results of
  Welters (see \cite{We}): The tangent map
  $\tau:\TT_{[C]}(\mathcal{Q})\to \TT_{[C-2h^2]}\JJ^2(X)$ is a map of
  vector spaces $\tau:\HH^0(C,N_{C/X})\to \HH^1(\Omega^2_X)^\vee$;
  using Serre duality and the isomorphism $N_{C/X}^\vee\isom
  N_{C/X}(-2)$, we get its dual map $\tau^\ast:\HH^1(\Omega^2_X) \to
  \HH^1(C,N_{C/X}(-1))$.  This map fits into a commutative diagram:

\[
\begin{array}{ccc}
\HH^0(X,\Oh_X(3)) & \by{R} & \HH^1(\Omega^2_X) \\
{r_C}\downarrow  & & \downarrow{\tau^\ast} \\
\HH^0(C,\Oh_C(3)) & \by{\partial} & \HH^1(C,N_{C/X}(-1)) \\
\end{array}
\]
where (see \cite{We} for details) $\partial$ is the coboundary
map in the long exact sequence of cohomologies associated to the short
exact sequence of normal bundles
$$0\to N_{C/X}(-1)\to N_{C/\bbP^4}(-1) \to \Oh_C(3) \to 0,$$  and as before
$R$ is the $\HH^0 \to \HH^1$ map of the sequence
$$0 \to \Omega^2_X \to \Omega^3_{\bbP^4}\tensor N_{X/\bbP^4} \to
 \Omega^3_X\tensor N_{X/\bbP^4} \to 0.$$

 We see that $R$, $r_C$ are surjections and the cokernel of $\partial$
 is $\HH^1(N_{C/\bbP^4}(-1))$. Since
 $N_{C/\bbP^4}=\Oh_C(2)^{\oplus{3}}$ ($C$ being a complete
 intersection of three quadrics in $\bbP^4$), the cokernel of
 $\partial$ in the diagram above has rank $3$. By Riemann-Roch,
 $h^0(N_{C/X})=\deg{N_{C/X}}+2(1-g(C))=2\cdot8-2\cdot 2\cdot
 2=8.$\comment{We note that $R$ is an isomorphism and $r_C$ is a
   surjection (from a $5$-dimensional vector space onto a
   $3$-dimensional vector space). We shall now show that
   $\dim{\im(\tau^\ast)}=2$. Since $\im{\tau^\ast}=\im{\partial}$, it
   is enough to prove that $\dim{\im(\partial)}=1$. Since $C$ is ACM,
   this implies that $C$ is a complete intersection in $\bbP^4$ and
   hence $N_{C/\bbP^4}$ is a sum of line bundles. Since
   $h^0(I_{C/X}(1))\neq 0$, this implies that $C\subset\bbP^3$, hence
   $N_{C/\bbP^4}$ already has $\Oh_C(1)$ as a summand. Now
   $\deg{C}=c_2(E)=2$, so $C\subset\bbP^3$ has to be a complete
   intersection of type $(1,2)$. This implies that
   $N_{C/\bbP^4}=\Oh_{\bbP^4}(1)^{\oplus{2}}\oplus\Oh_{\bbP^4}(2)$. From
   the associated cohomology sequence, one gets a surjection
   $\HH^1(C,N_{C/X}(-2)) \onto \HH^1(\Oh_C(-1))^{\oplus{2}}\isom
   \bbC^2$. Since $\HH^1(N_{C/X}(-2))\isom\HH^0(C,N_{C/X})$ is
   $4$-dimensional as seen before, we get $\im(\partial)$ is
   $2$-dimensional.} This implies that the map $\alpha$ has
 $3$-dimensional fibres and a $5$-dimensional image. But $\alpha$
 factors via $p:\mathcal{Q}\to\mathcal{M}$ and the fibres of $p$ are
 $\bbP(\HH^0(X,E))\isom\bbP^3$. This implies that the map
 $c_2:\mathcal{M}\to\JJ^2(X)$ is quasi-finite and \'etale. Thus we are
 done.
\end{proof}

\subsection{Stable Bundles with $(c_1,c_2)=(1,3)$}\label{1,3}
\comment{We get a $4\times 4$ matrix with entries $(a,b,c,d,e,f)=
  (l_1,l_2,l_3,c_1,c_2,c_3)$.}  

\begin{lemma}
Let $E$ be a normalised, stable bundle of rank two on a smooth quartic
threefold $X$ with $(c_1,c_2)=(1,3)$. Then $E$ is ACM.
\end{lemma}

\begin{proof}
The proof is similar to the cubic case. Let $S$ be a very general
hyperplane section as before (so $S$ is a $K3$ surface). By the
Noether-Lefschetz theorem, we have
$\Pic(S)\isom\bbZ[\Oh_S(1)]$. Furthermore, by a result of Maruyama
\cite{maruyama}, $E_{|S}$ is stable. Hence $h^0(E_{|S}(-1))=0$. Recall
Riemann-Roch for a quartic surface:
$$\chi(E_{|S})=2c_1^2-c_2+4.$$ Using this we get $h^0(E_{|S})\geq
3$.

Take any non-zero section $s$ of $E$. As before, let $\Gamma$ denote
the intersection of $S$ with the zero locus of $s$; hence $\Gamma$ is
a zero-dimensional subscheme of length three. From the short exact
sequence
\begin{equation}\label{eq4}
0\longrightarrow \Oh_S \stackrel{\cdot s\vert_S}{\longrightarrow}
E_{|S} \longrightarrow  I_{\Gamma/S}(1) \longrightarrow 0
\end{equation}
we see that $h^0(I_{\Gamma/S}(1))\geq 2$. From the short exact
sequence
$$
0\to I_{\Gamma/S}(1) \to \Oh_S(1) \to \Oh_{\Gamma}(1) \to
0
$$
we get a four-term exact sequence
$$0\to \HH^0(I_{\Gamma/S}(1)) \to \HH^0(\Oh_S(1))\isom\bbC^4 \to
\HH^0(\Oh_{\Gamma}(1))\isom\bbC^3 \to\HH^1(I_{\Gamma/S}(1))\to 0.$$
Since $\Oh_S(1)$ is very ample, sections of this line bundle separate
points and tangents. Hence $$\dim\im[\HH^0(\Oh_S(1))\to
\HH^0(\Oh_\Gamma(1))]\geq 2\, .$$ This implies that
$h^0(I_{\Gamma/S}(1))=2$ and hence $h^1(I_{\Gamma/S}(1))=1$. The long
exact sequence of cohomologies associated with the exact sequence in
Eq. \eqref{eq4} yields:
$$0\to \HH^1(E_{|S}) \to \HH^1(I_{\Gamma/S}(1)) \to \HH^2(\Oh_S) \to
\HH^2(E_{|S}) \to \cdots. $$ We have $h^2(E_{|S})=h^0(E_{|S}^\vee)=0$,
$h^2(\Oh_S)=1$ and from above $h^1(I_{\Gamma/S}(1))=1$. Hence
$h^1(E_{|S})=0$.

The proof that $h^1(E_{|S}(1))=0$ is similar: using Riemann-Roch, we
see that $h^0(E_{|S}(1))\geq 11$ and so
$h^0(I_{\Gamma/S}(2))=h^0(E_{|S}(1))-h^0(\Oh_S(1))\geq 7$. In the
four-term exact sequence
$$0\to \HH^0(I_{\Gamma/S}(2)) \to \HH^0(\Oh_S(2)) \to
\HH^0(\Oh_{\Gamma}(2)) \to\HH^1(I_{\Gamma/S}(2))\to 0,$$ we have
$h^0(\Oh_S(2))=10$ and $h^0(\Oh_{\Gamma}(2))=3$. Note that
$\Gamma\subset L\isom \bbP^1$ because $h^0(I_{\Gamma/S}(1))=2$. This
proves that the map $\HH^0(\Oh_S(2)) \to \HH^0(\Oh_{\Gamma}(2))$ is
surjective. Thus $h^1(I_{\Gamma/S}(2))=0$ and hence $h^1(E_{|S}(1))=0$.  
\comment{as $\HH^0(I_{L/\bbP^3}(2))\subsetneq\HH^0(\Oh_{\bbP^3}(2))$} 

Since $h^2(E_{|S})=0$, this implies that the Castelnuovo-Mumford
regularity of $E_{|S}$ is $2$. Hence $h^1(E_{|S}(n))=0$ for $n\geq
1$. By Serre duality, $h^1(E_{|S}(m))=0$ for $m\leq -2$. Since
$h^1(E_{|S})=h^1(E_{|S}(-1)=0$, this means that $E_{|S}$ is ACM and
thus we are done as before.
\end{proof}

For $X$ general, one can check that the moduli space
$\mathcal{M}_X(2;1,3)$ is smooth by using the following resolution
whose proof is similar to the cubic case and hence we omit it.

\begin{prop}
Any normalised stable bundle $E$ of rank two with $(c_1,c_2)=(1,3)$ on
a smooth quartic threefold is ACM and admits the following minimal
resolution:
\begin{equation}
0 \to \Oh_{\bbP^4}(-1)\oplus \Oh_{\bbP^4}(-3)^{\oplus{3}}
\to \Oh_{\bbP^4}(-2) \oplus \Oh_{\bbP^4}^{\oplus{3}} \to E
\to 0.
\end{equation}
\end{prop}

\comment{
\begin{proof}
Same as in the cubic case.
\end{proof}
}
\begin{thm}
The moduli space $\mathcal{M}_X(2;1,3)$ is isomorphic to $F(X)$, the
Fano variety of lines in $X$. Hence, for $X$ general, this moduli space  is 
a smooth irreducible curve.
\end{thm}

\begin{proof}
The proof of the isomorphism is similar to the proof of Theorem
\ref{firstthm}. When $X$ is general,
$\HH^1(N_{C/X})=\HH^2(X,\sEnd{E})=0$ .  Thus $F(X)$ is also smooth. By
Riemann-Roch, it now follows that $\dim{F(X)}=1$.  Irreducibility is
well known (see \cite{Collino}).
\end{proof}

\subsection{Bundles with $(c_1,c_2)=(0,2)$}
\comment{This is the case of conics.} By reasoning as above, it is
easy to see that such bundles are ACM and admit the following minimal
resolution:\comment{The $4\times 4$ matrix with entries
  $(a,b,c,d,e,f)=(c_1,c_2,q_1,q_2,l_1,l_2)$.}
\begin{equation}
0 \to \Oh_{\bbP^4}(-4)\oplus \Oh_{\bbP^4}(-3)^{\oplus{2}}\oplus
\Oh_{\bbP^4}(-2) \to \Oh_{\bbP^4}
\oplus\Oh_{\bbP^4}(-1)^{\oplus{2}}\oplus\Oh_{\bbP^4}(-2) \to E \to 0.
\end{equation}

\begin{thm}
  There is an isomorphism $$\mathcal{M}_X(2;0,2)\isom F_2(X),$$ where
  $F_2(X)$ is the Hilbert scheme of plane conics in $X$. Hence, for
  $X$ general, the moduli space $\mathcal{M}_X(2;0,2)$ is a smooth
  irreducible surface.
\end{thm}

\begin{proof}
>From the resolution of $E$, it is clear that $E$ has unique section
whose zero locus is a plane conic $C$. So the map
$$\mathcal{M}_X(2;0,2) \to F_2(X)$$ is just the one which takes
$[E]\mapsto [C]$. Just as in the cubic case, here too
$\Ext^1(I_{C/X},\Oh_X)\isom \bbC$. Thus this map is an
isomorphism. The vanishing of $\HH^2(X,\sEnd(E))$ follows from results
of \cite{MPR1}. It follows then that $F_2(X)$ is smooth
i.e. $h^1(N_{C/X})=0$. By Riemann-Roch, we have
$h^0(N_{C/X})-h^1(N_{C/X})=\deg{N}+2(1-g)$. Since $N_{C/X}\isom
E_{|C}$, we have $$\deg{N_{C/X}}=c_1(E)c_2(E)=0\, .$$ Finally, the
genus of $C$ is $0$. Thus we have that $F_2(X)$ is of dimension
$2$. Irreducibility of $F_2(X)$ has been proved in \cite{CMW},
smoothness follows from Lemma \ref{ismooth}.
\end{proof}

\subsection{Bundles with $(c_1,c_2)=(-1,1)$}
Recall that a rank two vector bundle $A$ with first Chern class $0$ or
$-1$, is {\it unstable} if and only if $h^0(A)=0$ (see \cite{OSS}). It
is easy to check using Grothendieck-Riemann-Roch that normalised,
indecomposable ACM bundles of rank two with $c_1=-1$ have a unique
section whose zero locus is a line i.e., $c_2=1$. Hence these bundles
are unstable. In fact, one can also check that any such bundle is of
the form $G(1)$ where for $E$ as in section \ref{1,3}, $G$ is defined
by the exact sequence
$$ 0 \to G \to \Oh_X(-2)\oplus\Oh_X^{\oplus{3}} \to E
\to 0.$$ 

Geometrically speaking, if we start with a line $\ell\subset X$, then
$G(1)$ is the rank two bundle which corresponds to $\ell$ via the
Serre construction. This can be seen as follows. Let $l_1, l_2, l_3$
be linear forms on $\bbP^4$ whose zero locus is $\ell$. Since
$\ell\subset X$, the quartic polynomial $f$ defining $X$ is of the form 
$$f=\sum_{i=1}^3  l_ic_i, $$ 
where the $c_i$'s are cubic polynomials. Let $\bbP^2\subset\bbP^4$ be
the plane defined by $l_1=0=l_2$. Then $\bbP^2\cap X$ is a reducible
curve $\ell\cup C$ where $C$ is the cubic curve given by
$l_1=l_2=c_3=0$. $G(1)$ and $E$ are the vector bundles which via
Serre's construction correspond to $\ell$ and $C$ respectively.

\section{An Application}

\subsection{Lagrangian Fibrations and complete integrable systems}
The results of this section are inspired by the observations of
Beauville in \cite{Beau2}. Let $X$ be a smooth cubic or quartic
threefold. Let $S\in|K_X^{-1}|$ be a general member so that $S$ is a
smooth $K3$ surface. By \cite{maruyama}, $E_{|S}$ is stable since
$\dim{X}>\rank{E}$. Let $\mathcal{M}_X$ (respectively,
$\mathcal{M}_S$) be any of the moduli spaces of stable vector bundles
of rank two on $X$ (respectively, $S$) discussed above and
$r:\mathcal{M}_X \to \mathcal{M}_S$ be the restriction map. Mukai
\cite{Mukai} has proved that $\mathcal{M}_S$ is smooth and carries a
symplectic structure. We have an exact sequence
$$0\to (\sEnd{E})\tensor{K}_X \to \sEnd{E} \to \sEnd{E}_{|S} \to 0.$$
Since $\HH^1(X,\sEnd{E})\isom \HH^2(X,\sEnd{E})=0$, on taking
cohomology we get an exact sequence
$$0 \to \HH^1(X,\sEnd{E}) \by{r_\ast} \HH^1(X,\sEnd{E}_{|S})
\by{r^\ast} \HH^1(X,\sEnd{E})^\vee \to 0,$$ where $r_\ast$ is the
tangent map to $r$ at $[E]$ and $r^\ast$ its transpose with respect to
the symplectic form. Therefore $r_\ast$ is injective and its image is
a maximal isotropic subspace. As a result, we have the following

\begin{prop}
$\mathcal{M}_X$ embeds as a Lagrangian submanifold of $\mathcal{M}_S$.
\end{prop}

As mentioned in \cite{Beau2}, if we fix an $S$ as above and consider
the linear system of hypersurfaces containing $S$, we get a Lagrangian
fibration ($=$ complete integrable system) as follows: Consider the
fibration $\mathcal{M}_S \to \Pi$ where $\Pi$ is the linear system of
hypersurfaces in $\bbP^4$ containing a given $X$. The fibre of this
map is $\mathcal{M}_X$. 

\comment{
\subsection{Lagrangian submanifolds and complete integrable systems}
\subsubsection{The cubic threefold}
Let $S\in|K_X^{-1}|$ be a general member so that $S=X\cap{Q}$ where
$Q\subset\bbP^4$ is a smooth quadric. By \cite{maruyama}, $E_{|S}$ is
stable since $\dim{X}>\rank{E}$. Then we have an exact sequence
$$0\to \sEnd{E}\tensor{K}_X \to \sEnd{E} \to \sEnd{E}_{|S} \to 0.$$
Since $\HH^1(X,\sEnd{E})\isom \HH^2(X,\sEnd{E})=0$, the corresponding
long exact sequence of cohomologies yields
$$0 \to \HH^1(X,\sEnd{E}) \to\HH^1(X,\sEnd{E}_{|S}) \to
\HH^1(X,\sEnd{E}) \to 0.$$

Let $\mathcal{M}_S$ be the moduli space of stable rank $2$ vector
bundles with $(c_1,c_2)=(1,2)$ and $r:\mathcal{M}_X \to \mathcal{M}_S$
be the restriction map. Then at any point $[E]\in\mathcal{M}_X$, we
have $T_{[E]}\mathcal{M}_X=\HH^1(X,\sEnd{E})$. The above sequence
implies the following

\begin{lemma}
$\mathcal{M}_X$ embeds as a Lagrangian submanifold of $\mathcal{M}_S$.
\end{lemma}

\subsubsection{The quartic threefold}
\begin{remark}
Notice that the map which sends $[E]\mapsto [E(-1)]$ induces an
isomorphism $\mathcal{M}_X(2;2,8) \isom \mathcal{M}_X(2;0,4)$. For
$S\in|K_X|$, we have fibration $\mathcal{M}_X(2;0.4) \to \Pi$ where
$\Pi$ is the linear system of quartic hypersurfaces in $\bbP^4$
containing a given $X$. The fibre of this map is
$\mathcal{M}_X(2;0,4)$ and this yields a complete integrable system
(see Beauville). As remarked in {\it op. cit.}, the moduli space
$\mathcal{M}_S(2;0,4)$ is a very interesting space and has been
considered by Kieran O'Grady.
\end{remark}
}
\section{Conclusions}
The results proved here give us a fairly complete understanding of ACM
vector bundles of rank two on general smooth hypersurfaces in
$\bbP^4$. As mentioned in the introduction, in the case when $X$ is a
general smooth quintic threefold in $\bbP^4$, there are only finitely
many isomorphism classes of rank two indecomposable ACM bundles on
$X$. In other words, there are only finitely many ways in which a
general smooth quintic can be realised as the Pfaffian of a minimal
skew-symmetric matrix of size $2k\times 2k$ for $k>1$. It would be
interesting to compute this number and understand its geometric
significance. In the case when the matrix has linear entries, this
number is ``an instance of the holomorphic Casson invariant'' (see
\cite{Beau}).

\section{Acknowledgements}
The third named author would like to thank N.~Mohan Kumar and
A.~Prabhakar Rao for educating him about various aspects and issues of
vector bundles and also for initiating him into the program of
studying ACM bundles on hypersurfaces through collaboration. The
questions addressed in this note first arose in discussions during
that collaboration.

\end{document}